# Equivalence of Weighted DT-Moduli of (Co)convex Functions

Malik Saad Al-Muhja[1,2,*]   Habibulla Akhadkulov[1,*]   Nazihah Ahmad[1,*,†]

**Abstract.**

The paper present new definitions for weighted DT moduli. Similarly, we a general outcome in an equivalence of moduli of smoothness are obtained. It is known that, any $r \in \mathbb{N}_o$, $0 < p \leq \infty$, $1 \leq \eta \leq r$ and $\phi(x) = \sqrt{1-x^2}$, the inequalities $\omega^{\phi}_{i+1,r}(f^{(r)}, \|\theta_{\mathcal{N}}\|)_{w_{\alpha,\beta,p}} \sim \omega^{\phi}_{i,r+1}(f^{(r+1)}, \|\theta_{\mathcal{N}}\|)_{w_{\alpha,\beta,p}}$ and $\omega^{\phi}_{i+\eta}(f, \|\theta_{\mathcal{N}}\|)_{\alpha,\beta,p} \sim \|\theta_{\mathcal{N}}\|^{-\eta} \omega^{\phi}_{i,2\eta}(f^{(2\eta)}, \|\theta_{\mathcal{N}}\|)_{\alpha+\eta,\beta+\eta,p}$ are valid.

## 1. Introduction

Hierarchy foundations of the moduli of smoothness are begun modern with the work of Ditzian and Totik (1987), (see [6]), and Kopotun (2006-2019), (see [8, 9, 10, 11, 12, 14, 15, 16, 18]). Ditzian and Totik established better continuous moduli of the function in a normed space, then Kopotun contributed to properties of various moduli of smoothness like univariate piecewise polynomial functions (splines) [16]. He has a significant impact on the hierarchy between moduli of smoothness for the past 14 years and we are affected by his contribution in the $k$th symmetric difference (see, [9, proof of Lemma 4.1]). Let $\Delta_h^k(f,x)$ be the $k$th symmetric difference of $f$ is given [6] by

$$\Delta_h^k(f,x) = \begin{cases} \sum_{i=0}^{k} \binom{k}{i}(-1)^{k-i} f\left(x + \left(\frac{2i-k}{2}\right)h\right), & x \pm \frac{kh}{2} \in [-1,1] \\ 0, & \text{otherwise.} \end{cases}$$

The space $L_p([-1,1])$, $0 < p < \infty$, denotes the space of all measurable functions $f$ on $[-1,1]$, [15] such that

$$\|f\|_{L_p[-1,1]} = \begin{cases} \left(\int_{-1}^{1} |f(x)|^p dx\right)^{\frac{1}{p}} < \infty, & 0 < p < \infty \\ \underset{x \in [-1,1]}{\text{esssup}} |f(x)|, & p = \infty. \end{cases}$$

---

[1] Department of Mathematics and statistics, School of Quantitative Sciences, College of Arts and Sciences, Universiti Utara Malaysia, 06010 Sintok, Kedah, Malaysia
[2] Department of Mathematics and Computer Application, College of Sciences, University of Al-Muthanna, Samawa 66001, Iraq

† **Corresponding author:** Nazihah Ahmad

* E-mail addresses: malik@mu.edu.iq ; dr.al-muhja@hotmail.com (M. Al-Muhja); habibulla@uum.edu.my (H. Ahadkulov); nazihah@uum.edu.my (N. Ahmad)

Let $\|\cdot\|_p = \|\cdot\|_{L_p[-1,1]}$, $0 < p \leq \infty$ and $\phi(x) = \sqrt{1-x^2}$. Then, Ditzian-Totik modulus of smoothness (DTMS) of a function $f \in L_p[-1,1]$, is defined [5] by

$$\omega_{k,r}^\phi(f,t)_p = \sup_{0 < h \leq t} \left\|\phi^r \Delta_{h\phi}^k(f,x)\right\|_p, \qquad k, r \in \mathbb{N}_o.$$

Also, the $k$th usually modulus of smoothness of $f \in L_p[-1,1]$ is defined [6] by

$$\omega_k(f,\delta,[-1,1])_p = \sup_{0 < h \leq \delta} \left\|\Delta_h^k(f,x)\right\|_p, \qquad \delta > 0, p \leq \infty.$$

Denote by $AC_{loc}(-1,1)$ and $AC[-1,1]$ the set of functions whose are locally absolutely continuous on $(-1,1)$ and absolutely continuous on $[-1,1]$ respectively. Now, we will need to accept the following:

**Definition 1.1 [18]** Let $w_{\alpha,\beta}(x) = (1+x)^\alpha (1-x)^\beta$ be the (classical) Jacobi weight, and let

$$\alpha, \beta \in J_p = \begin{cases} \left(\frac{-1}{p}, \infty\right), & \text{if } p < \infty, \\ [0, \infty), & \text{if } p = \infty. \end{cases}$$

Define

$$\mathbb{L}_p^{\alpha,\beta} = \left\{f : [-1,1] \to \mathbb{R} \mid \|w_{\alpha,\beta} f\|_p < \infty, \text{ and } 0 < p < \infty\right\},$$

$$\mathbb{L}_{p,r}^{\alpha,\beta} = \left\{f : [-1,1] \to \mathbb{R} \mid f^{(r-1)} \in AC_{loc}(-1,1), 1 \leq p \leq \infty \text{ and } \|w_{\alpha,\beta} f^{(r)}\|_p < \infty\right\}$$

and for convenience denote $\mathbb{L}_{p,0}^{\alpha,\beta} = \mathbb{L}_p^{\alpha,\beta}$.

Let $f \in \mathbb{L}_{p,r}^{\alpha,\beta}$, we write $\|\cdot\|_{w_{\alpha,\beta},p}$. If $r = 0$, we denoted $\|\cdot\|_{\alpha,\beta,p}$.

**Definition 1.2 [20]** A subset $X$ of $\mathbb{R}^n$ is convex set if $[x,y] \subseteq X$, whenever $x, y \in X$. Equivalently, $X$ is convex if

$$(1-\lambda)x + \lambda y \in X, \text{ for all } x, y \in X \text{ and } \lambda \in [0,1].$$

The function $f$ is called convex of $X$ if

$$f\big((1-\lambda)x + \lambda y\big) \leq (1-\lambda)f(x) + \lambda f(y), \text{ for all } x, y \in X \text{ and } \lambda \in [0,1].$$

**Definition 1.3 [7]** Let $Y_s = \{y_i\}_{i=1}^s$, $s \in \mathbb{N}$ be a partition of $[-1,1]$, that is, a collection of $s$ fixed points $y_i$ such that

$$y_{s+1} = -1 < y_s < \cdots < y_1 < 1 = y_o$$

and let $\Delta^{(2)}(Y_s)$ be the set of continuous functions on $[-1,1]$ that are convex downwards on the segment $[y_{i+1}, y_i]$ if $i$ is even and convex upwards on the same segment if $i$ is odd. The functions from $\Delta^{(2)}(Y_s)$ are called coconvex.

**Definition 1.4 [21]** The partition $\tilde{T}_\eta = \{t_j\}_{j=0}^{\eta}$, where

$$t_j = t_{j,\eta} = \begin{cases} -\cos\left(\dfrac{j\Pi}{\eta}\right), & \text{if } 0 \leq j \leq \eta, \\ -1, & \text{if } j < 0, \end{cases}$$

and $t_j$'s as the knots of Chebyshev partition.

**Definition 1.5 [17]** A function $f$ is said to be $k$-monotone, $k \geq 1$ on $[a,b]$, if and only if for all choices of $k+1$ distinct points $x_o, \ldots, x_k$ in $[a,b]$ the inequality $f[x_o, \ldots, x_k] \geq 0$, holds, where

$$f[x_o, \ldots, x_k] = \sum_{j=0}^{k} \frac{f(x_j)}{\theta'(x_j)}, \quad \theta(x_j) = \prod_{j=0}^{k}(x - x_j)$$

denotes the $k$th divided difference of $f$ at $x_o, \ldots, x_k$.

Now, we present the most important, Kopotun's methods and some further developments of his contribution in the $k$th symmetric difference. He said that (see, [16]): "$A$ is equivalent to $B$ and write $A \sim B$ if $c^{-1}A \leq B \leq cA$ such that $c$ positive constant". Let us recall:

**FIRST:** a piecewise polynomial $s$ on Chebyshev partition of $[-1,1]$, (see, [12]) then

$$\omega_{k+\eta}^{\phi}(s,t)_p \leq ct^\eta \omega_{k,\eta}^{\phi}(s^{(\eta)}, t)_p, \quad 0 < p < 1, t > 0,$$

and

$$\omega_{k-\eta,\eta}^{\phi}(s^{(\eta)}, n^{-1})_p \sim \omega_k^{\phi}(s, n^{-1})_p.$$

***Second:*** in 2007, Kopotun dedicated attention to the computation of several results on the equivalence of moduli of smoothness are obtained (see [16]), for example,

$$n^\eta \omega_{k-\eta}^{\phi}(s^{(\eta)}, n^{-1})_p \sim \omega_k(s, n^{-1})_p, \quad 1 \leq p \leq \infty, 1 \leq \eta \leq \min\{k, m+1\}. \tag{1}$$

***Third:*** the $k$-monotone functions were of major interest to Kopotun [14]. He has examined in 2009 the equivalence

$$\omega_k(f,\delta)_p \leq Ac_\delta(k,q,p)\|f\|_p,$$

where $f$ is satisfied (1), $q < p$ and

$$c_\delta(k,q,p) = \begin{cases} \delta^{\frac{2}{q}-\frac{2}{p}}, & \text{if } k \geq 2 \\ \delta^{\frac{2}{q}-\frac{2}{p}}, & \text{if } k = 1 \text{ and } p < 2q \\ \left(\delta\sqrt{|\ln(\delta)|}\right)^{\frac{1}{q}}, & \text{if } k = 1 \text{ and } p = 2q \\ \delta^{\frac{1}{q}}, & \text{if } k = 1 \text{ and } p > 2q. \end{cases}$$

The first to deal with development moduli of smoothness were Kopotun, Leviatan and Shevshuk [8]. They were interested with discuss various properties of the new modulus of smoothness

$$\omega_{k,r}^{\phi}(f^{(r)}, t)_p = \sup_{0 < h \leq t} \|\mathcal{W}_{kh}^r(.)\Delta_{h\phi}^k(f^{(r)},.)\|_p,$$

where

$$\mathcal{W}_{kh}^r(x) = \begin{cases} \left(\left(1-x-\delta\frac{\phi(x)}{2}\right)\left(1+x-\delta\frac{\phi(x)}{2}\right)\right)^{\frac{1}{2}}, & \text{if } 1 \pm x - \delta\frac{\phi(x)}{2} \in [-1,1] \\ 0, & \text{otherwise.} \end{cases}$$

However, they showed contribution in the $k$th symmetric difference of his modulus by K-functional which show immediately in [9, proof of Lemma 4.1].

The following result was proven by different method of modulus of smoothness in [11].

**Theorem 1.6** Let $k, n \in \mathbb{N}, r \in \mathbb{N}_o, A > 0, 0 < p \leq \infty, \left(\alpha + \frac{r}{2}\right), \left(\beta + \frac{r}{2}\right) \in J_p$, and let $0 < t \leq \varrho n^{-1}$, where $\varrho$ is some positive constant that depends only on $\alpha, \beta, k$ and $q$. Then, for any $p_n \in \pi_n$,

$$\omega_{k,r}^\phi\left(p_n^{(r)}, t\right)_{\alpha,\beta,p} \sim \Psi_{k,r}^\phi\left(p_n^{(r)}, t\right)_{\alpha,\beta,p} \sim \Omega_{k,r}^\phi\left(p_n^{(r)}, A, t\right)_{\alpha,\beta,p} \sim t^k \left\|w_{\alpha,\beta}\phi^r p_n^{(k+r)}\right\|_p,$$

where

$$\Psi_{k,r}^\phi\left(p_n^{(r)}, t\right)_{\alpha,\beta,p} = \sup_{0 \leq h \leq t} \left\|w_{\alpha,\beta}\phi^r \Delta_{h\phi}^k\left(p_n^{(r)}, x\right)\right\|_p,$$

$$\Omega_{k,r}^\phi\left(p_n^{(r)}, A, t\right)_{\alpha,\beta,p} = \sup_{0 \leq h \leq t} \left\|w_{\alpha,\beta}\phi^r \Delta_{h\phi}^k\left(p_n^{(r)}, x; \mathfrak{T}_{A,h}\right)\right\|_{L_p(\mathfrak{T}_{A,h})}$$

and the equivalence constants depend only on $k, r, \alpha, \beta, A$ and $q$.

**Definition 1.7 [10]** For $r \in \mathbb{N}_o$ and $0 < p \leq \infty$, denote $\mathbb{B}_p^0(w_{\alpha,\beta}) = \mathbb{L}_p^{\alpha,\beta}$ and

$$\mathbb{B}_p^{(r)}(w_{\alpha,\beta}) = \left\{f \mid f^{(r-1)} \in AC_{loc}(-1,1) \text{ and } \phi^r f^{(r)} \in \mathbb{L}_p^{\alpha,\beta}\right\}, \quad r \geq 1.$$

In 2018, Kopotun et al. ([10, Lemmas 2.2, 2.3]) proposed a function $f \in \mathbb{B}_p^{(r)}(w_{\alpha,\beta})$ and $\left(\frac{r}{2} + \alpha\right) \geq 0$, $\left(\frac{r}{2} + \beta\right) \geq 0$. Then,

$$\omega_{k,r}^\phi\left(f^{(r)}, t\right)_{\alpha,\beta,p} \leq c \left\|w_{\alpha,\beta}\phi^r f^{(r)}\right\|_p, \quad t > 0,$$

and

$$\lim_{t \to 0^+} \omega_{k,r}^\phi\left(f^{(r)}, t\right)_{\alpha,\beta,p} = 0.$$

The paper structure is as follows: In Section 2, our necessary results are stated. Then, our main results for weighted DT moduli appear in Section 3.

## 2. Notations and Further Results

In this section, we will present linear space for functions of Lebesgue Stieltjes integrable-i. Firstlly, let us recall the definition of the Lebesgue Stieltjes integrable-i, given in [3].

**Definition 2.1** Let $\mathbb{D}$ be measurable set, $f: \mathbb{D} \to \mathbb{R}$ be a bounded function, and $\mathcal{L}_i: \mathbb{D} \to \mathbb{R}$ be nondecreasing function for $i \in \Lambda$. For a Lebesgue partition P of $\mathbb{D}$, put $\underline{LS}(f, P, \underline{\mathcal{L}}) = \sum_{j=1}^n \prod_{i \in \Lambda} m_j \mathcal{L}_i\left(\mu(\mathbb{D}_j)\right)$ and $\overline{LS}(f, P, \underline{\mathcal{L}}) = \sum_{j=1}^n \prod_{i \in \Lambda} M_j \mathcal{L}_i\left(\mu(\mathbb{D}_j)\right)$ where $\mu$ is ameasure function of $\mathbb{D}$, $m_j = \inf\{f(x): x \in \mathbb{D}_j\}$, $M_j = \sup\{f(x): x \in \mathbb{D}_j\}$ and $\underline{\mathcal{L}} = \mathcal{L}_1, \mathcal{L}_2, \ldots$. Also, $\mathcal{L}_i(x_j) - \mathcal{L}_i(x_{j-1}) > 0$, $\underline{LS}(f, P, \underline{\mathcal{L}}) \leq \overline{LS}(f, P, \underline{\mathcal{L}})$, $\prod_{i \in \Lambda} \underline{\int}_i^{\mathbb{D}} f \underline{d\mathcal{L}} = \sup\{\underline{LS}(f, \underline{\mathcal{L}})\}$ and $\prod_{i \in \Lambda} \overline{\int}_i^{\mathbb{D}} f \underline{d\mathcal{L}} = \inf\{\overline{LS}(f, \underline{\mathcal{L}})\}$ where $\underline{LS}(f, \underline{\mathcal{L}}) = \{\underline{LS}(f, P, \underline{\mathcal{L}}): P \text{ part of set } \mathbb{D}\}$ and $\overline{LS}(f, \underline{\mathcal{L}}) = \{\overline{LS}(f, P, \underline{\mathcal{L}}): P \text{ part of set } \mathbb{D}\}$. If $\prod_{i \in \Lambda} \underline{\int}_i^{\mathbb{D}} f \underline{d\mathcal{L}} = \prod_{i \in \Lambda} \overline{\int}_i^{\mathbb{D}} f \underline{d\mathcal{L}}$ where $\underline{d\mathcal{L}} = d\mathcal{L}_1 \times d\mathcal{L}_2 \times \ldots$. Then $f$ is integral $\int_i$ according to $\mathcal{L}_i$ for $i \in \Lambda$.

**Lemma 2.2 [2]** If $f$ is a function of Lebesgue Stieltjes integral-i, then $vf$ is a function of Lebesgue Stieltjes integral-i, where $v > 0$ is real number, and

$$\prod_{i \in \Lambda} \int_i^{\mathbb{D}} vf \underline{d\mathcal{L}} = v \prod_{i \in \Lambda} \int_i^{\mathbb{D}} f \underline{d\mathcal{L}},$$

holds.

**Lemma 2.3 [2]** If the functions $f_1, f_2$ are integrable on the set $\mathbb{D}$ according to $\mathcal{L}_i$, for $i \in \Lambda$, then $f_1 + f_2$ is the function of integrable according to $\mathcal{L}_i$, for $i \in \Lambda$, such that

$$\prod_{i \in \Lambda} \int_i^{\mathbb{D}} (f_1 + f_2) \underline{d\mathcal{L}} = \prod_{i \in \Lambda} \int_i^{\mathbb{D}} f_1 \underline{d\mathcal{L}} + \prod_{i \in \Lambda} \int_i^{\mathbb{D}} f_2 \underline{d\mathcal{L}}.$$

**Definition 2.4 [4]** A domain $\mathbb{D}$ of convex polynomial $p_n$ of $\Delta^{(2)}$ is a subset of $X$ and $X \subseteq \mathbb{R}$, satisfying the following properties:

1) $\mathbb{D} \in \mathcal{K}^N$, where
$$\mathcal{K}^N = \{\mathbb{D}: \mathbb{D} \text{ is a compact subset of } X\}$$
is the class of all domain of convex polynomial,
2) there is the point $t \in X/\mathbb{D}$, such that
$$|p_n(t)| > \sup\{|p_n(x)|: x \in \mathbb{D}\}, \text{ and}$$
3) there is the function $f$ of $\Delta^{(2)}$, such that
$$\|f - p_n\| \leq \frac{c}{n^2} \omega_{2,2}^\phi \left(f'', \frac{1}{2}\right).$$

**Definition 2.5 [4]** A domain $\mathbb{D}$ of coconvex polynomial $p_n$ of $\Delta^{(2)}(Y_s)$ is a subset of $X$ and $X \subseteq \mathbb{R}$, satisfying the following properties:

1) $\mathbb{D} \in \mathcal{K}^N(Y_s)$, where
$$\mathcal{K}^N(Y_s) = \left\{\begin{matrix}\mathbb{D}: \mathbb{D} \text{ is a compact subset of } X, \\ \text{and } p_n \text{ changes convexity at } \mathbb{D}\end{matrix}\right\}$$
is the class of all domain of coconvex polynomial,
2) $y_i$'s are inflection points, such that
$$|p_n(y_i)| \leq \frac{1}{2}, i = 1, \ldots, s, \text{ and}$$
3) there is the function $f$ of $\Delta^{(2)}(Y_s)$, such that
$$\|f - p_n\| \leq \frac{c}{n^2} \omega_{k,2}^\phi \left(f'', \frac{1}{n}\right).$$

From Definitions 2.1, 2.4 and 2.5, if the function $f$ is convex, then $\mathbb{D}$ is domain of (co)convex function of $f$.

**Remark 2.6 [2]** Let $\mathcal{I}_f$ be denote the class of all functions of integrable $f$ that satisfying Definition 2.1, i.e.,

$$\mathcal{I}_f = \{f : f \text{ is integrable function according to } \mathcal{L}_i, i \in \Lambda\}$$

$$= \left\{ f : \prod_{i \in \Lambda} \underline{\int_i}^{\mathbb{D}} f \, d\underline{\mathcal{L}} = \prod_{i \in \Lambda} \overline{\int_i}^{\mathbb{D}} f \, d\underline{\mathcal{L}} \right\}.$$

**Remark 2.7 [13]** Let $x_i \in \left[\frac{x_i + x^\#}{2}, \frac{x_i + x_\star}{2}\right] \subseteq \theta_{\mathcal{N}}$, then we denote

$$x^\# = x_{j(i)+1}, \quad x_\star = x_{j(i)-2}$$

where

$$\theta_{\mathcal{N}} = \theta_{\mathcal{N}}[-1,1] = \{x_i\}_{i=0}^{\mathcal{N}} = \{-1 = x_\circ \leq \cdots \leq x_{\mathcal{N}-1} \leq x_{\mathcal{N}} = 1\}$$

and

$$\|\theta_{\mathcal{N}}\| = \max_{0 \leq i \leq \mathcal{N}-1} \{x_{i+1} - x_i\}$$

the length of the largest interval in that partition.

Next definitions are an immediate summary of [2] and [1].

**Definition 2.8** For $r \in \mathbb{N}_o$, the weighted DTMS in $\mathbb{L}_p^{\alpha,\beta} \cap \mathcal{I}_f$, we define

$$\Delta_h^i(f, x) = \begin{cases} \prod_{i \in \Lambda} \int_i^{\mathbb{D}} f \, d\underline{\mathcal{L}} & ; \text{ if } f \in \mathcal{I}_f \\ 0 & ; \quad \text{otherwise.} \end{cases} \quad (2)$$

By virtue of (2) and Definition 1.1, we define

$$\omega_{i,r}^\phi\left(f^{(r)}, \|\theta_{\mathcal{N}}\|, [-1,1]\right)_{w_{\alpha,\beta}, p} = \sup \left\{ \left\| w_{\alpha,\beta} \phi^r \Delta_{h\phi}^i(f^{(r)}, x) \right\|_p, 0 < h \leq \|\theta_{\mathcal{N}}\| \right\},$$

where $\|\theta_{\mathcal{N}}\| < 2(i^{-1})$, $\mathcal{N} \geq 2$.

**Definition 2.9** For $\alpha, \beta \in J_p$, $r \in \mathbb{N}_o$ and $0 < p \leq \infty$, we denote

$$\Phi^{p,r}(w_{\alpha,\beta}) = \left\{ f : f \in \mathbb{L}_{p,r}^{\alpha,\beta} \cap \mathcal{I}_f \text{ and } \omega_{i,r}^\phi\left(f^{(r)}, \|\theta_{\mathcal{N}}\|, [-1,1]\right)_{w_{\alpha,\beta}, p} < \infty \right\},$$

and $\Phi^{p,0}(w_{\alpha,\beta}) = \Phi^p(w_{\alpha,\beta})$.

Main contribution to this work, we focus to applications of results were obtained in last papers (see [2, Theorems 3.1, 3.3] and [1, Theorem 2.11]). More precisely, using some results of this paper like Definitions 2.8 and 2.9, we outcomes of direct estimates are obtaining.

A set of all piecewise polynomial approximation $\mathbb{S}(\tilde{T}_\eta, r+2)$ of order $r+2$, with the knots of Chebyshev partition $\tilde{T}_\eta$.

**Theorem 2.10 [2]** For $r \in \mathbb{N}_o$, $\alpha, \beta \in J_p$, there is a constant $c = c(r, \alpha, \beta, p)$ such that if $f \in \Delta^{(2)} \cap \mathbb{L}_{p,r}^{\alpha,\beta}$, there, a number $\mathcal{N} = \mathcal{N}\left(f, \omega_{1,r}^{\phi}(f^{(r)}, \|\theta_{\mathcal{N}}\|, I)_{w_{\alpha,\beta},p}\right)$ for $n \geq \mathcal{N}$ and $\mathcal{S} \in \mathbb{S}(\tilde{T}_\eta, r+2) \cap \Delta^{(2)} \cap \mathbb{L}_{p,r}^{\alpha,\beta}$, such that

$$\|f^{(r)} - \mathcal{S}^{(r)}\|_{w_{\alpha,\beta},p} \leq c_{r,\alpha,\beta,p,\omega_{1,r}^{\phi}} \min\left\{\omega_{i,r}^{\phi}(f^{(r)}, \|\theta_{\mathcal{N}}\|, I_\alpha)_{w_{\alpha,\beta},p}, \omega_{i,r}^{\phi}(f^{(r)}, \|\theta_{\mathcal{N}}\|, I_\beta)_{w_{\alpha,\beta},p}\right\}$$

where

$$\Delta_{h\phi,\alpha}^{i}(f^{(r)}, x) = \int_{1}^{\mathbb{D}} \int_{2}^{\mathbb{D}} \cdots \int_{i}^{\mathbb{D}} \cdots f^{(r)}(x)\, d\mathcal{L}_{1t,\alpha} d\mathcal{L}_{2t,\alpha} \cdots d\mathcal{L}_{it,\alpha} \cdots = \prod_{i \in \Lambda} \int_{i}^{\mathbb{D}} f^{(r)}\, \underline{d\mathcal{L}_{t\phi,\alpha}},$$

$$\Delta_{h\phi,\beta}^{i}(f^{(r)}, x) = \int_{1}^{\mathbb{D}} \int_{2}^{\mathbb{D}} \cdots \int_{i}^{\mathbb{D}} \cdots f^{(r)}(x)\, d\mathcal{L}_{1t,\beta} d\mathcal{L}_{2t,\beta} \cdots d\mathcal{L}_{it,\beta} \cdots = \prod_{i \in \Lambda} \int_{i}^{\mathbb{D}} f^{(r)}\, \underline{d\mathcal{L}_{t\phi,\beta}}.$$

Moreover, if $r, \alpha, \beta = 0$, then

$$\|f - \mathcal{S}\|_p \leq c\left(\omega_{1}^{\phi}\right) \omega_{i}^{\phi}(f, \|\theta_{\mathcal{N}}\|, I)_p.$$

In particular,

$$\|f^{(r)} - \mathcal{S}^{(r)}\|_{w_{\alpha,\beta},p} \leq c_r\, \omega_{1,r}^{\phi}(f^{(r)}, \|\theta_{\mathcal{N}}\|, I)_{w_{\alpha,\beta},p}.$$

**Theorem 2.11 [2]** Let $\Delta^k$ be the space of all $k$-monotone functions. If $f \in \Delta^k \cap \mathbb{L}_{p,r}^{\alpha,\beta}$ is such that $f^{(r)}(x) = p_n^{(r)}(x)$, where $p_n \in \pi_n \cap \Delta^k$, $N \geq k \geq 2$ and $s \in \mathbb{S}(\tilde{T}_\eta, r+2) \cap \Delta^k \cap \mathbb{L}_{p,r}^{\alpha,\beta}$. Then

$$\|f - s\|_{w_{\alpha,\beta},p} \leq c(f, p, k, \alpha, \beta, x_\star, x^\#)\omega_{i,r}^{\phi}(f, \|\theta_{\mathcal{N}}\|, I)_{w_{\alpha,\beta},p}.$$

In particular, if $f$ is a convex function and $p_n$ is a convex polynomial or piecewise convex polynomial, then

$$\|f - s\|_{w_{\alpha,\beta},p} \leq c_k \omega_{i,r}^{\phi}(f, \|\theta_{\mathcal{N}}\|, I)_{w_{\alpha,\beta},p}.$$

**Definition 2.12 [1]** For $\alpha, \beta \in J_p$ and $f \in \mathcal{I}_f$, we set

$$\mathbb{E}_n(f, w_{\alpha,\beta})_{\alpha,\beta,p} = \mathbb{E}_n(f)_{\alpha,\beta,p} = \inf\{\|f - p_n\|_{\alpha,\beta,p},\ p_n \in \pi_n \cap \mathcal{I}_f,\ f \in \Delta^{(2)}(Y_s) \cap \Phi^p(w_{\alpha,\beta})\}$$

and

$$\mathcal{E}_n^{(2)}(f, w_{\alpha,\beta}, Y_s)_p = \inf\{\|f - p_n\|_{\alpha,\beta,p},\ p_n \in \pi_n \cap \Delta^{(2)}(Y_s) \cap \mathcal{I}_f,\ f \in \Delta^{(2)}(Y_s) \cap \Phi^p(w_{\alpha,\beta})\}$$

respectively, denote the degree of best unconstrained and (co)convex polynomial approximation of $f$.

**Theorem 2.13 [1]** Let $\sigma, m, n \in \mathbb{N}$, $\sigma \neq 4$, $s \in \mathbb{N}_o$ and $\alpha, \beta \in J_p$. If $f \in \Delta^{(2)}(Y_s) \cap \Phi^p(w_{\alpha,\beta})$, then

$$\sup\left\{n^\sigma \mathcal{E}_n^{(2)}(f, w_{\alpha,\beta}, Y_s)_p : n \geq m\right\} \leq c \sup\{n^\sigma \mathbb{E}_n(f)_{\alpha,\beta,p} : n \in \mathbb{N}\}.$$

In particular, suppose that $Y_s \in \mathbb{Y}_s$ and $s \geq 1$. Then

$$\mathcal{E}_n^{(2)}(f, w_{\alpha,\beta}, Y_s)_p \leq cn^{-\sigma}\omega_{i,r}^{\phi}(f^{(r)}, \|\theta_{\mathcal{N}}\|, I)_{\alpha,\beta,p}, n \geq \|\theta_{\mathcal{N}}\|.$$

**Remark 2.14 [1]** If $f$ in $\mathcal{I}_f$ is a function of Lebesgue Stieltjes integral-i, and $f$ is differentiable function, therefore,

$$f' = \frac{df}{dx} = \frac{d}{dx}\left(\int_0^x \frac{df(u)}{d\mathcal{G}_{1,\mu,\mathbb{D}_o}} d\mathcal{G}_{1,\mu,\mathbb{D}_o}\right)$$

$$= \frac{d}{dx}\left(\int_0^x \int_0^x \frac{d^2 f(u)}{d\mathcal{G}_{1,\mu,\mathbb{D}_o} \times d\mathcal{G}_{2,\mu,\mathbb{D}_o}} d\mathcal{G}_{1,\mu,\mathbb{D}_o} \times d\mathcal{G}_{2,\mu,\mathbb{D}_o}\right)$$

$$= \frac{d}{dx}\left(\int_0^x \int_0^x \cdots \int_0^x \cdots \frac{d^i f(u)}{d\mathcal{G}_{1,\mu,\mathbb{D}_o} \times d\mathcal{G}_{2,\mu,\mathbb{D}_o} \times \ldots \times d\mathcal{G}_{i,\mu,\mathbb{D}_o} \times \ldots} d\mathcal{G}_{1,\mu,\mathbb{D}_o} \times d\mathcal{G}_{2,\mu,\mathbb{D}_o} \times \ldots \times d\mathcal{G}_{i,\mu,\mathbb{D}_o} \times \ldots\right)$$

$$= \frac{d}{dx}\left(\prod_{i \in \Lambda} \int_i^{I_x} f^{(i)}(u)\, \underline{d\mathcal{G}_{\mu,\mathbb{D}_o}}\right), x \in I_x = [0, x] \subseteq \mathbb{D}_o, u \in \mathbb{D}_o \text{ and } \mathcal{G}_{\mu,\mathbb{D}_o} = \mathcal{G}(\mu(\mathbb{D}_o))$$

$$= \frac{d}{dx}\left(\prod_{i \in \Lambda} \int_i^{I_x} f^{(i)}\, \underline{d\mathcal{G}_{\mu,\mathbb{D}_o}}\right)$$

$$= \prod_{i \in \Lambda} \int_i^{I_x} f_x^{(i+1)}\, \underline{d\mathcal{G}_{\mu,\mathbb{D}_o}}.$$

$$f_x^{(i+1)} = \frac{d^i}{d\mathcal{G}_{i,\mu,\mathbb{D}_o}^i} f' = \frac{d}{dx}\left(\frac{d^i f}{d\mathcal{G}_{i,\mu,\mathbb{D}_o}^i}\right) = \frac{d^{i+1} f_x}{dx \times d\mathcal{G}_{i,\mu,\mathbb{D}_o}^i}.$$

By Definition 2.9, we show the following result.

**Lemma 2.15** We have

$$\Phi^{p,r+1}(w_{\alpha,\beta}) = \Phi^{p,r}\left(w_{\alpha+\frac{1}{2},\beta+\frac{1}{2}}\right).$$

Proof. Firstly, suppose $1 \leq p < \infty$, and $w_{\alpha,\beta}(x) = (1+x)^{\alpha}(1-x)^{\beta}$.

Let $f \in \Phi^{p,r+1}(w_{\alpha,\beta})$ and assume $f$ satisfy Definition 2.8. Next,

$$\|w_{\alpha,\beta}\phi^{r+1}\Delta_{h\phi}^i(f^{(r+1)}, x)\|_p = \left(\int_{-1}^1 |w_{\alpha,\beta}\phi^{r+1}\Delta_{h\phi}^i(f^{(r+1)}, x)|^p dx\right)^{\frac{1}{p}}, 0 < h \leq \|\theta_{\mathcal{N}}\|$$

$$= \left(\int_{-1}^1 \left|w_{\alpha,\beta}\phi^{r+1}\prod_{i \in \Lambda}\int_i^{\mathbb{D}} f^{(r+1)}\, \underline{d\mathcal{L}_{\phi}}\right|^p dx\right)^{\frac{1}{p}}.$$

Finally, from [2, proof of Lemma 3.2], [18] and Remark 2.14, then

□

$$\left\|w_{\alpha,\beta}\phi^{r+1}\Delta_{h\phi}^{i}(f^{(r+1)},x)\right\|_{p} = \left(\int_{-1}^{1}\left|w_{\alpha+\frac{1}{2},\beta+\frac{1}{2}}\phi^{r}\prod_{i\in\Lambda}\int_{i+1}^{\mathbb{D}}f^{(r)}\,\underline{d\mathcal{L}_{\phi}}\right|^{p}dx\right)^{\frac{1}{p}}$$

$$= \left\|w_{\alpha+\frac{1}{2},\beta+\frac{1}{2}}\phi^{r}\Delta_{h\phi}^{i+1}(f^{(r)},x)\right\|_{p}, 0 < h \leq \|\theta_{\mathcal{N}}\|.$$

**Remark 2.16** By virtue of Lemma 2.15, we see the following immediate consequence

$$\omega_{i,r+1}^{\phi}(f^{(r+1)},\|\theta_{\mathcal{N}}\|)_{w_{\alpha,\beta},p} = \omega_{i+1,r}^{\phi}(f^{(r)},\|\theta_{\mathcal{N}}\|)_{w_{\alpha+\frac{1}{2},\beta+\frac{1}{2}},p}.$$

**3. Main Results for Weighted DT Moduli**

In this section, we have the following main results.

**Theorem 3.1** Let $s, r \in \mathbb{N}_o$, $0 < p \leq \infty$ and $f \in \Delta^{(2)}(Y_s) \cap \Phi^{p,r}(w_{\alpha,\beta})$. Let $\mathbb{D}$ be defined in Definition 2.1, such that $|\mathbb{D}| \leq \delta_o$, for some $\delta_o \in \mathbb{R}^+$. Then,

$$\omega_{i+1,r}^{\phi}(f^{(r)},\|\theta_{\mathcal{N}}\|)_{w_{\alpha,\beta},p} \leq c(\delta_o)\omega_{i,r+1}^{\phi}(f^{(r+1)},\|\theta_{\mathcal{N}}\|)_{w_{\alpha,\beta},p}, \tag{3}$$

where the constant $c$ depend on $\delta_o$.

Proof. Suppose that $f \in \Delta^{(2)}(Y_s) \cap \Phi^{p,r}(w_{\alpha,\beta})$,

$$\oint_{i}^{\mathbb{D}_k \cap \mathbb{D}_j} = \prod_{i\in\Lambda}\int_{i}^{\mathbb{D}_k \cap \mathbb{D}_j} f\,\underline{d\mathcal{L}_{\phi}}$$

and

$$\oint_{i}^{\mathbb{D} \cap \mathbb{D}_j} = \prod_{i\in\Lambda}\int_{i}^{\mathbb{D} \cap \mathbb{D}_j} f\,\underline{d\mathcal{L}_{\phi}}.$$

We also, assume that $\mathbb{D}_j \subset \mathbb{D}$ such that

$$f(x) = \begin{cases} |\mathbb{D}| & \text{if } |\mathbb{D}| \leq \delta_o \\ \left(\oint_{i}^{\mathbb{D}_k \cap \mathbb{D}_j}\right) \to \left(\oint_{i}^{\mathbb{D} \cap \mathbb{D}_j}\right) & \text{if } \mathbb{D}_k, \mathbb{D}_j \text{ are sets of Lebeguse measurable,} \\ 0, & \text{otherwise.} \end{cases} \tag{4}$$

Then,

$$\left\|w_{\alpha,\beta}\phi^{r}\Delta_{h\phi}^{i+1}(f^{(r)},x)\right\|_{p} = \left(\int_{-1}^{1}\left|w_{\alpha,\beta}\phi^{r}\prod_{i\in\Lambda}\int_{i+1}^{\mathbb{D}}f^{(r)}\,\underline{d\mathcal{L}_{\phi}}\right|^{p}dx\right)^{\frac{1}{p}}$$

$$= \begin{cases} \left(\int_{-1}^{1}\left|w_{\alpha,\beta}\phi^r \prod_{i\in\Lambda}\int_{i+1}^{\mathbb{D}}|\mathbb{D}|\,\underline{d\mathcal{L}_\phi}\right|^p dx\right)^{\frac{1}{p}} = \mathrm{I}_o(x) & \text{if } |\mathbb{D}|\leq \delta_o,\, \delta_o\in\mathbb{R}^+, \\ \left(\int_{-1}^{1}\left|w_{\alpha,\beta}\phi^r \lim_{k\to\infty}\prod_{i\in\Lambda}\int_{i+1}^{\mathbb{D}_k\cap\mathbb{D}_j} f^{(r)}\,\underline{d\mathcal{L}_\phi}\right|^p dx\right)^{\frac{1}{p}} = \mathrm{I}_1(x) & \text{if } \mathbb{D}_k,\, \mathbb{D}_j \text{ are sets of Lebeguse measurable,} \\ 0, & \text{otherwise.} \end{cases}$$

Therefore, (4) implies that

$$\mathrm{I}_o(x) = \left(\int_{-1}^{1}\left|w_{\alpha,\beta}\phi^r \prod_{i\in\Lambda}\int_{i+1}^{\mathbb{D}}|\mathbb{D}|\,\underline{d\mathcal{L}_\phi}\right|^p dx\right)^{\frac{1}{p}} \leq \delta_o$$

for some $\delta_o \in \mathbb{R}^+$, while

$$\mathrm{I}_1(x) = \left(\int_{-1}^{1}\left|w_{\alpha,\beta}\phi^r \lim_{k\to\infty}\prod_{i\in\Lambda}\int_{i+1}^{\mathbb{D}_k\cap\mathbb{D}_j} f^{(r)}\,\underline{d\mathcal{L}_\phi}\right|^p dx\right)^{\frac{1}{p}}$$

$$= \left(\int_{-1}^{1}\left|w_{\alpha,\beta}\phi^r \lim_{k\to\infty} LS\left(f^{(r)}, \underline{\mathcal{L}_\phi\left(\mu(\mathbb{D}_k\cap\mathbb{D}_j)\right)}\right)\right|^p dx\right)^{\frac{1}{p}}$$

$$= \left(\int_{-1}^{1}\left|w_{\alpha,\beta}\phi^r LS\left(f^{(r)}, \underline{\mathcal{L}_\phi\left(\lim_{k\to\infty}\mu(\mathbb{D}_k\cap\mathbb{D}_j)\right)}\right)\right|^p dx\right)^{\frac{1}{p}}$$

$$= \left(\int_{-1}^{1}\left|w_{\alpha,\beta}\phi^r LS\left(f^{(r)}, \mathcal{L}_\phi\left(\mu\left(\underline{\left(\bigcup_{k=1}^{\infty}\mathbb{D}_k\right)\cap\mathbb{D}_j}\right)\right)\right)\right|^p dx\right)^{\frac{1}{p}}$$

$$= \left(\int_{-1}^{1}\left|w_{\alpha,\beta}\phi^r LS\left(f^{(r)}, \underline{\mathcal{L}_\phi\left(\mu(\mathbb{D}\cap\mathbb{D}_j)\right)}\right)\right|^p dx\right)^{\frac{1}{p}}$$

$$= \left(\int_{-1}^{1}\left|w_{\alpha,\beta}\phi^r \prod_{i\in\Lambda}\int_{i+1}^{\mathbb{D}\cap\mathbb{D}_j} f^{(r)}\,\underline{d\mathcal{L}_\phi}\right|^p dx\right)^{\frac{1}{p}}.$$

By Remark 2.14, we have

$$\mathrm{I}_1(x) \leq c\left(\int_{-1}^{1}\left|w_{\alpha,\beta}\phi^r \prod_{i\in\Lambda}\int_{i}^{\mathbb{D}\cap\mathbb{D}_j} f^{(r+1)}\,\underline{d\mathcal{L}_\phi}\right|^p dx\right)^{\frac{1}{p}}.$$

Taking supremum, we obtain (3). $\square$

This implies the following result is valid.

**Corollary 3.2** Let $s, r \in \mathbb{N}_o$, $0 < p \leq \infty$ and $f \in \Delta^{(2)}(Y_s) \cap \Phi^{p,r}(w_{\alpha,\beta})$. Let $\mathbb{D}$ and $\delta_o$ be defined in Theorem 3.1. Then,

$$\omega^{\phi}_{i+1,r}\left(f^{(r)}, \|\theta_{\mathcal{N}}\|\right)_{w_{\alpha,\beta},p} \leq c(\delta_o)\omega^{\phi}_{i+1,r}\left(f^{(r)}, \|\theta_{\mathcal{N}}\|\right)_{w_{\alpha+\frac{1}{2},\beta+\frac{1}{2}},p},$$

where the constant $c$ depend on $\delta_o$.

Proof. Clearly. □

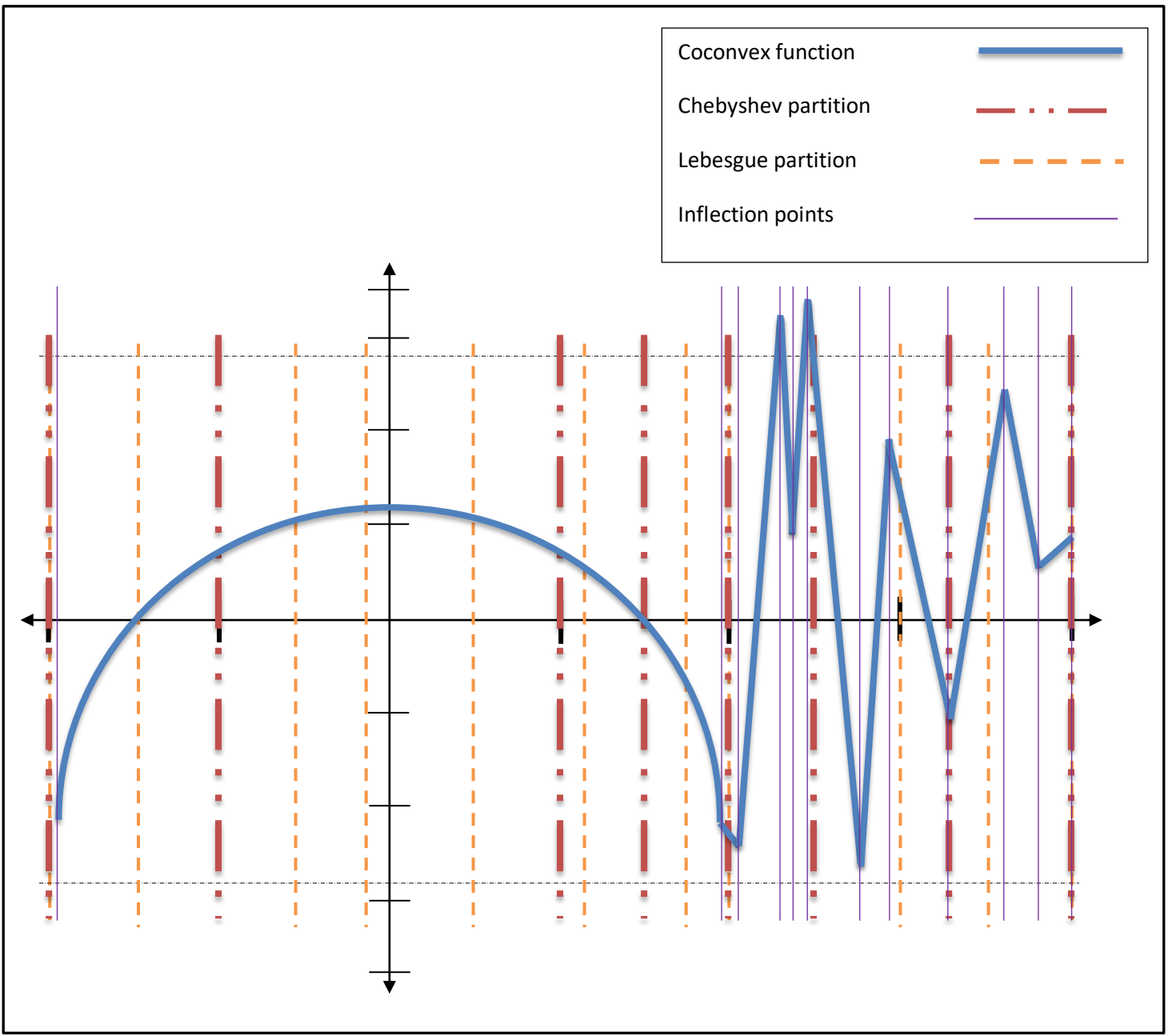

**Figure 1.** Show graph of partitions with interval $[-1,2]$ of the coconvex function.

**Theorem 3.3** Let $s, r \in \mathbb{N}_o$, $\alpha, \beta \in J_p$ and $0 < p \leq \infty$. Let P be a Lebesgue partition of $\mathbb{D}$, and $\tilde{T}_\eta$ be a Chebyshev partition with $P \cap \tilde{T}_\eta \neq \emptyset$, $1 \leq \eta \leq r$. If $f \in \Delta^{(2)}(Y_s) \cap \Phi^{p,r}(w_{\alpha,\beta})$, then there is a constant $c$ depend on $\eta$ and $J_{j,\eta}$ such that

$$\omega^\phi_{i+\eta}(f, \|\theta_\mathcal{N}\|)_{w_{\alpha,\beta},p} \leq c \|\theta_\mathcal{N}\|^{-\eta} \omega^\phi_{i,2\eta}(f^{(2\eta)}, \|\theta_\mathcal{N}\|)_{w_{\alpha+\eta,\beta+\eta},p}. \tag{5}$$

Proof. Recall that P is Lebesgue partition of $\mathbb{D}$, and $\tilde{T}_\eta$ is Chebyshev partition. Since $P \cap \tilde{T}_\eta \neq \emptyset$, virtue of [2, proof of Lemma 2.3], for $\varepsilon > 0$, then there is $P_\varepsilon$ Lebesgue partition of $\mathbb{D}$, which union with $\tilde{T}_\eta$ such that $P_\varepsilon \cup \tilde{T}_\eta = P$. We can create $J_{j,\eta} = [u_{j-(\eta+i)}, u_{j-(\eta+i)+1}]$ for some $y_i \in \cup_{j=0}^\eta J_{j,\eta}$ and $y_i$'s inflection points of $Y_s$, $s \in \mathbb{N}_o$, (see Figure 1). Next, if $f \in \Delta^{(2)}(Y_s) \cap \Phi^{p,r}(w_{\alpha,\beta})$, then,

$$\omega^\phi_{i+\eta}(f, \|\theta_\mathcal{N}\|)^p_{w_{\alpha,\beta},p} = \sup\left\{\left\|w_{\alpha,\beta} \phi^r \Delta^{i+\eta}_{h\phi}(f, x)\right\|^p_p, 0 < h \leq \|\theta_\mathcal{N}\|\right\}$$

$$\leq c \sup\left\{\sum_{j=0}^\eta \left\|w_{\alpha,\beta} \phi^r \Delta^{i+\eta}_{h\phi}(f - f^{(\eta)} + f^{(\eta)}, x)\right\|^p_{L_p(J_{j,\eta})}, 0 < h \leq \|\theta_\mathcal{N}\|\right\}.$$

By virtue of [19] and Theorem 2.13 or (see [1, proof of Theorem 2.11], implies

$$\|\theta_\mathcal{N}\|^\eta \omega^\phi_{i+\eta}(f, \|\theta_\mathcal{N}\|)^p_{w_{\alpha,\beta},p} \leq c \left(\sum_{j=0}^\eta \left(\int_{-1}^1 \left|w_{\alpha,\beta} \phi^r \left(\prod_{i \in \Lambda} \int_{i+\eta}^{J_{j,\eta}} (f - f^{(\eta)} + f^{(\eta)}) \, d\mathcal{L}_\phi\right)\right|^p dx\right)\right)$$

$$\leq c \left(\sum_{j=0}^\eta \left(\int_{-1}^1 \left|w_{\alpha,\beta} \phi^r \left(\prod_{i \in \Lambda} \int_{i+\eta}^{J_{j,\eta}} (f - f^{(\eta)}) \, d\mathcal{L}_\phi + \prod_{i \in \Lambda} \int_{i+\eta}^{J_{j,\eta}} f^{(\eta)} \, d\mathcal{L}_\phi\right)\right|^p dx\right)\right)$$

$$\leq c \sup\left\{\sum_{j=0}^\eta \left(\left\|w_{\alpha,\beta} \phi^r \Delta^{i+\eta}_{h\phi}(f - f^{(\eta)}, x)\right\|^p_{L_p(J_{j,\eta})} + \left\|w_{\alpha,\beta} \phi^r \Delta^{i+\eta}_{h\phi}(f^{(\eta)}, x)\right\|^p_{L_p(J_{j,\eta})}\right), 0 < h \leq \|\theta_\mathcal{N}\|\right\}$$

$$\leq c \left(\sup\left\{\sum_{j=0}^\eta \left\|w_{\alpha,\beta} \phi^r \Delta^{i+\eta}_{h\phi}(f - f^{(\eta)}, x)\right\|^p_{L_p(J_{j,\eta})}, 0 < h \leq \|\theta_\mathcal{N}\|\right\}\right.$$

$$\left.+ \sup\left\{\sum_{j=0}^\eta \left\|w_{\alpha,\beta} \phi^r \Delta^{i+\eta}_{h\phi}(f^{(\eta)}, x)\right\|^p_{L_p(J_{j,\eta})}, 0 < h \leq \|\theta_\mathcal{N}\|\right\}\right)$$

$$\leq c(\eta, J_{j,\eta}) \times \sup\left\{\sum_{j=0}^\eta \left\|w_{\alpha,\beta} \phi^r \Delta^{i+\eta}_{h\phi}(f^{(\eta)}, x)\right\|^p_{L_p(J_{j,\eta})}, 0 < h \leq \|\theta_\mathcal{N}\|\right\}$$

$$\leq c(\eta, J_{j,\eta}) \omega^\phi_{i+\eta,\eta}(f^{(\eta)}, \|\theta_\mathcal{N}\|)^p_{w_{\alpha,\beta},p}.$$

Now, by (3), implies (5) is proved. □

The following is an immediate of consequence of Theorems 3.1 and 3.3.

**Corollary 3.4** Let $s, r \in \mathbb{N}_o$, $\alpha, \beta \in J_p$, $0 < p \leq \infty$ and $f \in \Delta^{(2)}(Y_s) \cap \Phi^{p,r}(w_{\alpha,\beta})$. Let P be a Lebesgue partition of $\mathbb{D}$, and $\tilde{T}_\eta$ be a Chebyshev partition with $P \cap \tilde{T}_\eta \neq \emptyset$, $1 \leq \eta \leq r$. We have

$$\|w_{\alpha,\beta}\phi^\eta f^{(\eta)}\|_p \geq c(\eta, J_{j,\eta}) \begin{cases} \omega^\phi_{i+2\eta,i+\eta}(f^{(i+\eta)}, \|\theta_\mathcal{N}\|)_{w_{\alpha,\beta},p} & \text{if } |\mathbb{D}| \leq c(\eta, J_{j,\eta}), \\ \omega^\phi_{i,i+2\eta}(f^{(i+2\eta)}, \|\theta_\mathcal{N}\|)^p_{w_{\alpha+\frac{\eta}{2},\beta+\frac{\eta}{2}},p} & \text{if } |\mathbb{D}| > c(\eta, J_{j,\eta}). \end{cases} \quad (6)$$

Proof. Let $s, r \in \mathbb{N}_o$, $1 \leq \eta \leq r$, $\phi(x) = \sqrt{1-x^2}$ and $J_{j,\eta} = [u_{j-(\eta+i)}, u_{j-(\eta+i)+1}]$. Let P be a Lebesgue partition of $\mathbb{D}$, and $\tilde{T}_\eta$ be a Chebyshev partition. Assume the function $f \in \Delta^{(2)}(Y_s) \cap \Phi^{p,r}(w_{\alpha,\beta})$, and the constant $c$ depend on $\phi, r$ and $\eta$. Then,

$$c \times \|w_{\alpha,\beta}\phi^\eta f^{(\eta)}\|^p_p \geq \|w_{\alpha,\beta}\phi^\eta f^{(\eta)}\|^p_p$$

$$\geq \left\|w_{\alpha,\beta}\phi^\eta \left(\frac{\phi^r}{\phi^r}\right) f^{(\eta)}\right\|^p_p \geq c(\phi^r, \phi^\eta)^{-1} \|w_{\alpha,\beta}\phi^r f^{(\eta)}\|^p_p$$

$$\geq c(\phi^r, \phi^\eta)^{-1} \left\|w_{\alpha,\beta}\phi^r \prod_{\substack{i \in \Lambda \\ 1 \leq \eta \leq r}} \int_{i+\eta}^{J^x_{j,\eta}} f^{((i+\eta)+\eta)} \, d\mathcal{L}_\phi\right\|^p_p$$

$$\geq c(\phi^r, \phi^\eta)^{-1} \left\|w_{\alpha,\beta}\phi^r \prod_{\substack{i \in \Lambda \\ 1 \leq \eta \leq r}} \int_{i+\eta}^{J^x_{j,\eta}} f^{(i+2\eta)} \, d\mathcal{L}_\phi\right\|^p_p$$

$$\geq c(\eta, J_{j,\eta}) \sum_{j=0}^{\eta} \sup \left\|w_{\alpha,\beta}\phi^r \prod_{\substack{i \in \Lambda \\ 1 \leq \eta \leq r}} \int_{i+\eta}^{J^x_{j,\eta}} f^{(i+2\eta)} \, d\mathcal{L}_\phi\right\|^p_{L_p(J_{j,\eta})}$$

$$\geq c(\eta, J_{j,\eta}) \sup \left\{\sum_{j=0}^{\eta} \|w_{\alpha,\beta}\phi^r \Delta^{i+\eta}_{h\phi}(f^{(i+2\eta)}, x)\|^p_{L_p(J_{j,\eta})}, 0 < h \leq \|\theta_\mathcal{N}\|\right\}$$

$$\geq c(\eta, J_{j,\eta}) \omega^\phi_{i+\eta,i+2\eta}(f^{(i+2\eta)}, \|\theta_\mathcal{N}\|)^p_{w_{\alpha,\beta},p}.$$

Finally, by virtue of Theorems 3.1 and 3.3, we have

$$\|w_{\alpha,\beta}\phi^\eta f^{(\eta)}\|^p_p \geq c(\eta, J_{j,\eta}) \begin{cases} \omega^\phi_{i+2\eta,i+\eta}(f^{(i+\eta)}, \|\theta_\mathcal{N}\|)^p_{w_{\alpha,\beta},p} & \text{if } |\mathbb{D}| \leq c(\eta, J_{j,\eta}), \\ \omega^\phi_{i,i+2\eta}(f^{(i+2\eta)}, \|\theta_\mathcal{N}\|)^p_{w_{\alpha+\frac{\eta}{2},\beta+\frac{\eta}{2}},p} & \text{if } |\mathbb{D}| > c(\eta, J_{j,\eta}). \end{cases} \quad \square$$

## 4. Conclusions and Direct Estimates

The following theorem includes all outcomes of this paper. However, in this section, we state and provide equivalence corollaries for all outcomes in Section 3, that we provide.

**Theorem 4.1** Assume that $s, r \in \mathbb{N}_o$, $\alpha, \beta \in J_p$, $0 < p \leq \infty$ and $f \in \Delta^{(2)} \cap \Phi^{p,r}(w_{\alpha,\beta})$. If P is Lebesgue partition of $\mathbb{D}$, and $\tilde{T}_\eta$ is Chebyshev partition with $\mathbb{D} \cap \tilde{T}_\eta \neq \emptyset$. Then, for any constant $c$ may be depend on $\eta$ and $J_{j,\eta}$ and may be depend on $|\mathbb{D}| \leq \delta_o$, we have

$$\omega_{i+1,r}^\phi\big(f^{(r)}, \|\theta_\mathcal{N}\|\big)_{w_{\alpha,\beta},p} \sim c(\delta_o)\omega_{i,r+1}^\phi\big(f^{(r+1)}, \|\theta_\mathcal{N}\|\big)_{w_{\alpha,\beta},p} \sim$$

$$c(\delta_o) \times \omega_{i+1,r}^\phi\big(f^{(r)}, \|\theta_\mathcal{N}\|\big)_{w_{\alpha+\frac{1}{2},\beta+\frac{1}{2}},p} \sim \big\|w_{\alpha,\beta}\phi^\eta f^{(\eta)}\big\|_p \sim$$

$$c(\eta, J_{j,\eta})\left\{\omega_{i+2\eta,i+\eta}^\phi\big(f^{(i+\eta)}, \|\theta_\mathcal{N}\|\big)_{w_{\alpha,\beta},p} : \text{ if } |\mathbb{D}| \leq c(\eta, J_{j,\eta})\right\}$$

and

$$\|\theta_\mathcal{N}\|^\eta \times \omega_{i+\eta}^\phi(f, \|\theta_\mathcal{N}\|)_{w_{\alpha,\beta},p} \sim c(\eta, J_{j,\eta})\omega_{i,2\eta}^\phi\big(f^{(2\eta)}, \|\theta_\mathcal{N}\|\big)_{w_{\alpha+\eta,\beta+\eta},p} \sim$$

$$\big\|w_{\alpha,\beta}\phi^\eta f^{(\eta)}\big\|_p \sim c(\eta, J_{j,\eta})\left\{\omega_{i,i+2\eta}^\phi\big(f^{(i+2\eta)}, \|\theta_\mathcal{N}\|\big)_{w_{\alpha+\frac{\eta}{2},\beta+\frac{\eta}{2}},p}^p : \text{ if } |\mathbb{D}| > c(\eta, J_{j,\eta})\right\}.$$

**Corollary 4.2 ($s = 0$)** For $r \in \mathbb{N}_o$ and $\alpha, \beta \in J_p$, there is a constant $c$ may be depend on $r, \alpha, \beta, p, \omega_{1,r}^\phi$ and may be depend on $r, \alpha, \beta, p, \omega_{1,r}^\phi, \eta$ and $J_{j,\eta}$ such that $f \in \Delta^{(2)} \cap \Phi^{p,r}(w_{\alpha,\beta})$, $J_{j,\eta} = [u_{j-(\eta+i)}, u_{j-(\eta+i)+1}]$ and $1 \leq \eta \leq r$. Then,

$$\mathcal{E}_n^{(2)}\big(f, w_{\alpha,\beta}, Y_o\big)_p \leq c\|\theta_\mathcal{N}\|^\eta \omega_{i+\eta}^\phi(f, \|\theta_\mathcal{N}\|)_{w_{\alpha,\beta},p}$$

and

$$\mathcal{E}_n^{(2)}\big(f, w_{\alpha,\beta}, Y_o\big)_p \leq c(\eta, J_{j,\eta})\omega_{i,2\eta}^\phi\big(f^{(2\eta)}, \|\theta_\mathcal{N}\|\big)_{w_{\alpha+\eta,\beta+\eta},p}.$$

**Corollary 4.3 ($s \geq 1$)** Suppose that $Y_s \in \mathbb{Y}_s$, $\sigma, s, n \in \mathbb{N}$ and $\sigma \neq 4$. If $f \in \Delta^{(2)}(Y_s) \cap \Phi^{p,r}(w_{\alpha,\beta})$, then

$$\mathcal{E}_n^{(2)}\big(f, w_{\alpha,\beta}, Y_s\big)_p \leq c(\delta_o) \times n^{-\sigma}\omega_{i+1,r}^\phi\big(f^{(r)}, \|\theta_\mathcal{N}\|\big)_{w_{\alpha+\frac{1}{2},\beta+\frac{1}{2}},p}$$

and

$$\mathcal{E}_n^{(2)}\big(f, w_{\alpha,\beta}, Y_s\big)_p \leq c(\eta, J_{j,\eta}) \times n^{-\sigma}\omega_{i+2\eta,i+\eta}^\phi\big(f^{(i+\eta)}, \|\theta_\mathcal{N}\|\big)_{w_{\alpha,\beta},p}.$$

## 5. Data Availability

No data were used to support this study.

## 6. Conflict of Interests

The authors declare that there is no conflict of interests regarding the publication of this paper.

## 7. Acknowledgement

The first author is indebted to Prof. Eman Samir Bhaya (University of Babylon) for useful discussions of the subject. The first author is supported by University of Al-Muthanna while studying for his Ph. D. We would like to thank Universiti Utara Malaysia (UUM) for the financial support.